\documentclass[11pt]{article}
\usepackage{amssymb}

\usepackage{amsmath}
\usepackage{theorem}
\usepackage{here}
\usepackage[dvipdfmx]{color}

\newtheorem{thm}{Theorem}[section]
\newtheorem{lem}{Lemma}[section]
\newtheorem{cor}{Corollary}[section]
\newtheorem{prp}{Proposition}[section]
\theorembodyfont{\rmfamily}

\makeatletter

\@addtoreset{equation}{section}
\makeatother

\usepackage{comment} % 
\usepackage{bm}
\usepackage[pdftex]{graphicx}

\def\ta{{\tau}}

\def\al{{\alpha}}

\def\ga{{\gamma}}

\def\la{{\lambda}}

\def\om{{\omega}}
\def\th{{\theta}}

\def\ze{{\zeta}}

\def\bmu{{\text{\boldmath $\mu$}}}

\def\De{{\Delta}}

\def\Ga{{\Gamma}}

\def\x{{\text{\boldmath $x$}}}

\def\X{{\text{\boldmath $X$}}}
\def\Y{{\text{\boldmath $Y$}}}

\def\ph{{\hat p}}

\def\Tt{{\tilde T}}
\def\Ut{{\tilde U}}
\def\Vt{{\tilde V}}
\def\Wt{{\tilde W}}
\def\[{{\text{\boldmath $[$}}}
\def\]{{\text{\boldmath $]$}}}

\def\/{{\Bigr/\!\!}}

\def\1r{{\rm (1)}}
\def\2r{{\rm (2)}}
\def\3r{{\rm (3)}}
\def\4r{{\rm (4)}}
\def\5r{{\rm (5)}}

\def\non{{\nonumber}}
\begin{document}
\title{Bayesian Predictive Density Estimation for a Chi-squared Model Using Information from %
a Normal Observation with Unknown Mean and Variance}
\author{
Yasuyuki Hamura\footnote{Graduate School of Economics, University of Tokyo, 
7-3-1 Hongo, Bunkyo-ku, Tokyo 113-0033, JAPAN. 
JSPS Research Fellow.\newline{
E-Mail: yasu.stat@gmail.com}} \
and
Tatsuya Kubokawa\footnote{Faculty of Economics, University of Tokyo, 
7-3-1 Hongo, Bunkyo-ku, Tokyo 113-0033, JAPAN. \newline{
E-Mail: tatsuya@e.u-tokyo.ac.jp}}
}
%
\begin{comment}
%\\
{\it\normalsize  University of Tokyo}
%}
\end{comment}
%
\maketitle
\begin{abstract}
In this paper, we consider the problem of estimating the density function of a Chi-squared variable on the basis of observations of another Chi-squared variable and a normal variable under the Kullback-Leibler divergence. 
We assume that these variables have a common unknown scale parameter and that the mean of the normal variable is also unknown. 
We compare the risk functions of two Bayesian predictive densities: one with respect to a hierarchical shrinkage prior and the other based on a noninformative prior. 
The hierarchical Bayesian predictive density depends on the normal variable while the Bayesian predictive density based on the noninformative prior does not. 
Sufficient conditions for the former to dominate the latter are obtained. 
These predictive densities are compared by simulation. 
\par\vspace{4mm}
\end{abstract}

\section{Introduction}
\label{sec:introduction}
Suppose that $\X $ and $V$ are independently distributed according to the normal and Chi-squared distributions ${\rm{N}}_p ( \bmu , ( r_0 / \eta ) I_p )$ and $( r_0' / \eta ) \chi ^{2} ( n_1 )$ with densities 
\begin{align}
&p( \x | \bmu , \eta ) = {( \eta / r_0 )^{p / 2} \over (2 \pi )^{p / 2}} \exp \Big( - {\eta / r_0 \over 2} || \x - \bmu ||^2 \Big) \text{,} \quad \x \in \mathbb{R} ^p \text{,} \quad \text{and} \non \\
&p_1 (v | \eta ) = {(1 / 2)^{n_1 / 2} \over \Ga ( n_1 / 2)} v^{n_1 / 2 - 1} ( \eta / r_0' )^{n_1 / 2} \exp \Big( - {\eta / r_0' \over 2} v \Big) \text{,} \quad v \in (0, \infty ) \text{,} \non 
\end{align}
respectively, for known $p \in \mathbb{N} = \{ 1, 2, \dotsc \}$ and $r_0 , r_0' , n_1 > 0$ and unknown $\bmu \in \mathbb{R} ^p$ and $\eta \in (0, \infty )$. 
Suppose that for known $s_0' , n_2 > 0$, $W$ is an unobservable Chi-squared variable with distribution $( s_0' / \eta ) \chi ^{2} ( n_2 )$ which is independent of $( \X , V)$. 
We consider the problem of estimating the density of $W$, namely 
\begin{align}
p_2 (w | \eta ) = {(1 / 2)^{n_2 / 2} \over \Ga ( n_2 / 2)} w^{n_2 / 2 - 1} ( \eta / s_0' )^{n_2 / 2} \exp \Big( - {\eta / s_0' \over 2} w \Big) \text{,} \quad w \in (0, \infty ) \text{,} \non 
\end{align}
on the basis of the observation of $( \X , V)$ under the Kullback-Leibler loss. 
The risk function of a predictive density $\ph _2 ( \cdot ; \X , V)$ is 
\begin{align}
R(( \bmu , \eta ), \ph _2 ) &= E_{( \bmu , \eta )}^{( \X , V, W)} \Big[ \log {p_2 (W | \eta ) \over \ph _2 (W; \X , V)} \Big] \text{.} \non 
\end{align}
Such a situation arises, for example, if $\X _1 , \dots , \X _{N_1}$ and $\Y _1 , \dots , \Y _{N_2}$ are independently distributed as ${\rm{N}}_p ( \bmu , (1 / \eta ) I_p )$ and if we want to estimate the predictive density of $\sum_{i = 1}^{N_2} || \Y _i - \overline{\Y } ||^2$, where $\overline{\Y } = (1 / N_2 ) \sum_{i = 1}^{N_2} \Y _i$, based on the sufficient statistics $\overline{\X } = (1 / N_1 ) \sum_{i = 1}^{N_1} \X _i$ and $\sum_{i = 1}^{N_1} || \X _i - \overline{\X } ||^2$. 
On the other hand, since $r_0' / \eta $ and $n_1$ can be any positive real numbers, $V$ may be viewed as a gamma variable. 
Throughout the paper, however, we assume that $r_0 = r_0' = s_0' = 1$ for simplicity. 

For a prior $\pi ( \bmu , \eta )$ for the unknown parameters $( \bmu , \eta )$, the associated Bayesian predictive density $\ph _{2}^{( \pi )} ( \cdot ; \X , V)$ is given by 
\begin{align}
\ph _{2}^{( \pi )} (w; \x , v) &= E_{\pi }^{( \bmu , \eta ) | ( \X , V)} [ p_2 (w | \eta ) | ( \X , V) = ( \x , v) ] \non \\
&= {(1 / 2)^{n_2 / 2} \over \Ga ( n_2 / 2)} w^{n_2 / 2 - 1} E_{\pi }^{\eta | ( \X , V)} \Big[ \eta ^{n_2 / 2} \exp \Big( - {\eta \over 2} v \Big) \Big| ( \X , V) = ( \x , v) \Big] \text{.} \non 
\end{align}
The Jeffreys prior for the model where only $V$ is observed is $\pi _0 ( \bmu , \eta ) = \eta ^{- 1}$, %
which corresponds to the unbiased estimator $V / n_1$ of the variance $1 / \eta $ in the sense that $1 / E_{\pi _0} [ \eta | \X , V] = V / n_1$. 
As in Liang and Barron (2004), it can be shown that $\ph _{2}^{( \pi _0 )} ( \cdot ; \X , V)$ is uniformly optimal among the predictive densities which are equivariant with respect to the transformations of Section 2 of Stein (1964). 
In particular, for any $a_0 < n_1 / 2$, it improves upon $\ph _{2}^{( \pi _{a_0} )} ( \cdot ; \X , V)$ for $\pi _{a_0} ( \bmu , \eta ) = \eta ^{- a_0 - 1}$, which, when $a_0 = - p / 2$, coincides with the Jeffreys prior for the present model where both $\X $ and $V$ are observed. 
In this paper, as in %
Maruyama and Strawderman (2012), we consider the hierarchical shrinkage prior 
\begin{align}
\pi _{b, a} ( \bmu , \eta ) &= \int_{0}^{1} \pi _{b, a} ( \bmu , \ga , \eta ) d\ga \text{,} \label{eq:hierarchical_prior} 
\end{align}
where 
\begin{align}
\pi _{b, a} ( \bmu , \ga , \eta ) &= {\rm{N}}_p ( \bmu | \bm{0} _p , [ \{ (1 - \ga ) / \ga \} / \eta ] I_p ) (1 - \ga )^{b - 1} \ga ^{- a - 1} \eta ^{- a - 1} \non \\
&= {(1 - \ga )^{b - p / 2 - 1} \ga ^{p / 2 - a - 1} \eta ^{p / 2 - a - 1} \over (2 \pi )^{p / 2}} \exp \Big( - {\eta \over 2} {\ga \over 1 - \ga } || \bmu ||^2 \Big) \non 
\end{align}
for $b > 0$ and $a < p / 2$. 
We compare the two predictive densities $\ph _{2}^{( \pi _0 )} ( \cdot ; \X , V)$ and $\ph _{2}^{( \pi _{b, a} )} ( \cdot ; \X , V)$. 
In particular, in Section \ref{sec:dominance}, we obtain conditions under which $\ph _{2}^{( \pi _{b, a} )} ( \cdot ; \X , V)$ dominates $\ph _{2}^{( \pi _0 )} ( \cdot ; \X , V)$. 

An important feature of the problem is that the distribution of $\X $ depends on the unknown location parameter $\bmu $ while the distribution of $W$ does not depend on $\bmu $. 
As will be shown later, $\ph _{2}^{( \pi _0 )} ( \cdot ; \X , V)$ is a function only of $V$ but $\ph _{2}^{( \pi _{b, a} )} ( \cdot ; \X , V)$ does depend on $\X $. 
Thus, dominance of $\ph _{2}^{( \pi _{b, a} )} ( \cdot ; \X , V)$ over $\ph _{2}^{( \pi _0 )} ( \cdot ; \X , V)$ is analogous to the result of Stein (1964) that when estimating the variance $1 / \eta $ under the standardized squared error loss, the unbiased estimator $V / n_1$ can be improved upon by using additional information from $\X $. 

Although Stein (1964) considered a truncated estimator, it was shown by Brewster and Zidek (1974) that the unbiased estimator is dominated by a smooth generalized Bayes estimator also. 
Kubokawa (1994) showed that these improved estimators can be derived through the unified method of Integral Expression of Risk Difference (IERD). 
Maruyama (1998) gave a class of priors including that of Brewster and Zidek (1974) to improve on the unbiased estimator when the mean of the normal distribution is equal to zero. 
Related hierarchical priors have been shown to be useful in estimating location parameters in the presence of an unknown scale parameter (Maruyama and Strawderman (2005, 2020a, 2020b)). 

Bayesian predictive densities have been widely studied in the literature since Aitchison (1975) showed their superiority to plug-in predictive densities. 
Komaki (2001) proved for a normal model with unknown mean that the Bayesian predictive density against the uniform prior is dominated by that against a shrinkage prior as in estimation problems. 
Parallels between estimation and prediction were investigated by George, Liang and Xu (2006, 2012) and Brown, George and Xu (2008) in terms of minimaxity and admissibility. 
Kato (2009) and Boisbunon and Maruyama (2014) considered the case of unknown mean and variance. 
Prediction for a $2 \times 2$ Wishart model was considered by Komaki (2009). 
Prediction for a gamma model when the scale parameter is restricted to an interval was considered by L'Moudden, Marchand, Kortbi and Strawderman (2017).

\def\pit{{\tilde \pi}}

\section{Bayesian Predictive Densities}
\label{sec:predictive_densities}
In this section, the Bayesian predictive densities with respect to the priors $\pi _0$ and $\pi _{b, a}$ given in Section \ref{sec:introduction} are derived. 
The choice of the hyperparameter $b$ in $\pi _{b, a}$ is discussed. 

We first consider $\ph _{2}^{( \pi _0 )} ( \cdot ; \X , V)$ for the noninformative prior $\pi _0 ( \bmu , \eta ) = \eta ^{- 1}$. 

\begin{prp}
\label{prp:predictive_density_1} 
The Bayesian predictive density $\ph _{2}^{( \pi _0 )} ( \cdot ; \X , V)$ is given by 
\begin{align}
\ph _{2}^{( \pi _0 )} (w; \X , V) &= {1 \over B( n_1 / 2, n_2 / 2)} {V^{n_1 / 2} w^{n_2 / 2 - 1} \over (V + w)^{( n_1 + n_2 ) / 2}} \text{.} \non 
\end{align}
\end{prp}

\noindent
We note that this predictive density does not depend on $\X $. 
Moreover, it is identical to the predictive density with respect to the observation $V \sim (1 / \eta ) \chi ^2 ( n_1 )$ and the prior $\eta \sim \eta ^{- 1}$. 
Its superiority to the corresponding plug-in predictive density is discussed in Aitchison (1975). 

On the other hand, $\ph _{2}^{( \pi _{b, a} )} ( \cdot ; \X , V)$ actually depends on the normal variable $\X $. 

\begin{prp}
\label{prp:predictive_density_2} 
The Bayesian predictive density $\ph _{2}^{( \pi _{b, a} )} ( \cdot ; \X , V)$ for the hierarchical prior $\pi _{b, a}$ in (\ref{eq:hierarchical_prior}) is given by 
\begin{align}
\ph _{2}^{( \pi _{b , a} )} (w | \X , V) %
&= {w^{n_2 / 2 - 1} \over B( n_1 / 2 + p / 2 - a, n_2 / 2)} \frac{ \displaystyle \int_{0}^{1} {(1 - \ga )^{b - 1} \ga ^{p / 2 - a - 1} \over (V + w + \ga || \X ||^2 )^{( n_1 + n_2 ) / 2 + p / 2 - a}} d\ga }{ \displaystyle \int_{0}^{1} {(1 - \ga )^{b - 1} \ga ^{p / 2 - a - 1} \over (V + \ga || \X ||^2 )^{n_1 / 2 + p / 2 - a}} d\ga } \text{.} \non 
\end{align}
\end{prp}

\noindent
Because of the integrals in the above expression, the risk function of $\ph _{2}^{( \pi _{b, a} )} ( \cdot ; \X , V)$ is hard to evaluate in general. 

If we choose $b = n_1 / 2$, then the integral in the denominator can be simplified to 
\begin{align}
{B( n_1 / 2, p / 2 - a) \over V^{n_1 / 2} (V + || \X ||^2 )^{p / 2 - a}} \label{eq:denominator} 
\end{align}
by Lemma 2 of Boisbunon and Maruyama (2014). 
This choice corresponds to that in Section 2.1 of Maruyama and Strawderman (2005). 
On the other hand, in this case, the integral in the numerator becomes, by Lemma 2 of Boisbunon and Maruyama (2014), 
\begin{align}
{1 \over (V + w)^{( n_1 + n_2 ) / 2} (V + w + || \X ||^2 )^{p / 2 - a}} \int_{0}^{1} (1 - \ga )^{n_1 / 2 - 1} \ga ^{p / 2 - a - 1} \Big( 1 - {|| \X ||^2 \over V + w + || \X ||^2} \ga \Big) ^{n_2 / 2} d\ga \label{eq:numerator} 
\end{align}
and involves the hypergeometric function, which shows the greater complexity of the prediction problem. 
However, the above integral can be evaluated as in the proof of Lemma A2 of Boisbunon and Maruyama (2014), which is crucial for our proof of Theorem \ref{thm:BM} for general $n_2$. 

There is another case where we can analytically examine the risk function of $\ph _{2}^{( \pi _{b, a} )} ( \cdot ; \X , V)$. 
Suppose that $b = 1$. 
Then $\ph _{2}^{( \pi _{b, a} )} ( \cdot ; \X , V)$ becomes, by Lemma \ref{lem:incomplete_beta} in the Appendix, 
\begin{align}
\ph _{2}^{( \pi _{1, a} )} (w; \X , V) &= {w^{n_2 / 2 - 1} \over B( n_1 / 2 + p / 2 - a, n_2 / 2)} \frac{ \displaystyle \int_{0}^{1} {\ga ^{p / 2 - a - 1} \over (V + w + \ga || \X ||^2 )^{( n_1 + n_2 ) / 2 + p / 2 - a}} d\ga }{ \displaystyle \int_{0}^{1} {\ga ^{p / 2 - a - 1} \over (V + \ga || \X ||^2 )^{n_1 / 2 + p / 2 - a}} d\ga } \non \\
&= \ph _{2}^{( \pi _0 )} (w; \X , V) \frac{ \displaystyle \int_{0}^{|| \X ||^2 / (V + w + || \X ||^2 )} {\ga ^{p / 2 - a - 1} (1 - \ga )^{( n_1 + n_2 ) / 2 - 1} \over B(p / 2 - a, ( n_1 + n_2 ) / 2)} d\ga }{ \displaystyle \int_{0}^{|| \X ||^2 / (V + || \X ||^2 )} {\ga ^{p / 2 - a - 1} (1 - \ga )^{n_1 / 2 - 1} \over B(p / 2 - a, n_1 / 2)} d\ga } \text{.} \label{eq:predictive_density_2_1al} 
\end{align}
Therefore, 
\begin{align}
\lim_{|| \x ||^2 \to \infty } \ph _{2}^{( \pi _{1, a} )} ( \cdot ; \x , V) %
&= \ph _{2}^{( \pi _0 )} (w; \X , V) \text{,} \non 
\end{align}
which shows that we can apply the method of IERD of Kubokawa (1994). 
In order to prove Theorem \ref{thm:KK} given later, we use the expression (\ref{eq:predictive_density_2_1al}) and apply the argument of Kato (2009). 
Finally, %
it is interesting to note that for $a = p / 2 - 1$, the Bayesian predictive density $\ph _{2}^{( \pi _{1, a} )} ( \cdot ; \X , V)$ can be expressed in closed form as 
\begin{align}
\ph _{2}^{( \pi _{1, p / 2 - 1} )} ( \cdot ; \X , V) &= {n_1 + n_2 \over n_1 \Ga ( n_2 / 2)} {V^{n_1 / 2} w^{n_2 / 2 - 1} \over (V + w)^{( n_1 + n_2 ) / 2}} \frac{ \displaystyle 1 - \Big( {V + w \over V + w + || \X ||^2} \Big) ^{( n_1 + n_2 ) / 2} }{ \displaystyle 1 - \Big( {V \over V + || \X ||^2} \Big) ^{n_1 / 2} } \text{.} \label{eq:predictive_density_2_1al_special_case} 
\end{align}
That we can obtain this simple estimator is one of the important features of our prediction problem.

\section{Dominance Conditions%
}
\label{sec:dominance}
In this section, we provide sufficient conditions for $\ph _{2}^{( \pi _{b, a} )} ( \cdot ; \X , V)$ to dominate $\ph _{2}^{( \pi _0 )} ( \cdot ; \X , V)$ in the two cases $b = n_1 / 2$ and $b = 1$. 
In particular, conditions on the other hyperparameter $a$ are obtained. 

We first consider the case $b = n_1 / 2$. 
Let 
\begin{align}
( c_1 , c_2 ) &= \begin{cases} \displaystyle \Big( {\Ga ( n_1 / 2) \Ga (( n_1 + n_2 ) / 2 + p / 2 - a) \over \Ga (( n_1 + n_2 ) / 2) \Ga ( n_1 / 2 + p / 2 - a)} - 1, 1 \Big) \text{,} & \text{if $n_2 \le 2$} \text{,} \\ \displaystyle \Big( {p / 2 - a \over ( n_1 + n_2 ) / 2 - 1}, {n_2 \over 2} \Big) \text{,} & \text{if $n_2 > 2$} \text{.} \end{cases} \non 
\end{align}

\begin{thm}
\label{thm:BM} 
Suppose that $b = n_1 / 2$ and $a < p / 2$. 
If the inequality 
\begin{align}
&{p / 2 - a \over c_2} \Big\{ \psi \Big( {n_1 + n_2 \over 2} + {p \over 2} \Big) - \psi \Big( {n_1 \over 2} + {p \over 2} \Big) \Big\} \non \\
&\le \int_{0}^{1} (1 - \rho )^{( n_1 + n_2 ) / 2 + p / 2 - 1} {1 \over \rho } \Big\{ 1 - {1 \over (1 + c_1 \rho )^{( n_1 + n_2 ) / 2}} \Big\} d\rho \label{eq:condition_BM} 
\end{align}
is satisfied, then $R(( \bmu , \eta ), \ph _{2}^{( \pi _{b, a} )} ) \le R(( \bmu , \eta ), \ph _{2}^{( \pi _0 )} )$ for all $\bmu \in \mathbb{R} ^p$ and $\eta \in (0, \infty )$. 
Equality can hold only if $\bmu = \bm{0} _p$. 
\end{thm}

The integral appearing in the right-hand side of (\ref{eq:condition_BM}) is not a big problem. 
First, we can numerically calculate the integral since it does not involve the unknown parameters. 
Second, the integral can actually be evaluated analytically to obtain simpler sufficient conditions. 

\begin{cor}
\label{cor:D0} 
Assume that $b = n_1 / 2$ and $a < p / 2$. 
\begin{itemize}
\item[{\rm{(i)}}]
If 
\begin{align}
&\psi \Big( {n_1 + n_2 \over 2} + {p \over 2} \Big) - \psi \Big( {n_1 \over 2} + {p \over 2} \Big) \non \\
&\le {c_2 \over p / 2 - a} {n_1 + n_2 + p + 2 \over n_1 + n_2 + p} \Big[ 1 - {1 \over \{ 1 + 2 c_1 / ( n_1 + n_2 + p + 2) \} ^{( n_1 + n_2 ) / 2}} \Big] \text{,} \non 
\end{align}
then $\ph _{2}^{( \pi _{b, a} )} ( \cdot ; \X , V)$ dominates $\ph _{2}^{( \pi _0 )} ( \cdot ; \X , V)$. 
\item[{\rm{(ii)}}]
Suppose that either $n_2 \le 2$ and 
\begin{align}
\psi \Big( {n_1 + n_2 \over 2} + {p \over 2} \Big) - \psi \Big( {n_1 \over 2} + {p \over 2} \Big) < {( n_1 + n_2 ) c_2 \over n_1 + n_2 + p} \psi \Big( {n_1 + n_2 \over 2} \Big) - \psi \Big( {n_1 \over 2} \Big) \non 
\end{align}
or $n_2 > 2$ and 
\begin{align}
\psi \Big( {n_1 + n_2 \over 2} + {p \over 2} \Big) - \psi \Big( {n_1 \over 2} + {p \over 2} \Big) < {( n_1 + n_2 ) c_2 \over n_1 + n_2 + p} {2 \over n_1 + n_2 - 2} \text{.} \non 
\end{align}
Then $\ph _{2}^{( \pi _{b, a} )} ( \cdot ; \X , V)$ dominates $\ph _{2}^{( \pi _0 )} ( \cdot ; \X , V)$ for any $0 \le a < p / 2$ sufficiently close to $p / 2$. 
\end{itemize}
\end{cor}

When $n_2 = 2$, condition (\ref{eq:condition_BM}) is actually necessary and sufficient for $\ph _{2}^{( \pi _{n_1 / 2, a} )} ( \cdot ; \X , V)$ to dominate $\ph _{2}^{( \pi _0 )} ( \cdot ; \X , V)$. 

\begin{cor}
\label{cor:2n_2} 
Assume that $b = n_1 / 2$, $a < p / 2$, and $n_2 = 2$. 
\begin{itemize}
\item[{\rm{(i)}}]
$\ph _{2}^{( \pi _{b, a} )} ( \cdot ; \X , V)$ dominates $\ph _{2}^{( \pi _0 )} ( \cdot ; \X , V)$ if and only if 
\begin{align}
&{p / 2 - a \over n_1 / 2 + p / 2} \le \int_{0}^{1} (1 - \rho )^{n_1 / 2 + p / 2} {1 \over \rho } \Big( 1 - {1 \over [1 + \{ (p / 2 - a) / ( n_1 / 2) \} \rho ]^{n_1 / 2 + 1}} \Big) d\rho \text{.} \label{eq:c2n_2-1} 
\end{align}
\item[{\rm{(ii)}}]
When $n_1 = 2$, $\ph _{2}^{( \pi _{b, a} )} ( \cdot ; \X , V)$ dominates $\ph _{2}^{( \pi _0 )} ( \cdot ; \X , V)$ if and only if $0 \le a < p / 2$. 
\end{itemize}
\end{cor}

Next we consider the case of $b = 1$. 

\begin{thm}
\label{thm:KK} 
Assume that $b = 1$, $0 \le a < p / 2$, and $n_1 > 2$. 
Then $R(( \bmu , \eta ), \ph _{2}^{( \pi _{b, a} )} ) \le R(( \bmu , \eta ), \ph _{2}^{( \pi _0 )} )$ for all $\bmu \in \mathbb{R} ^p$ and $\eta \in (0, \infty )$. 
Equality holds if and only if $\bmu = \bm{0} _p$ and $a = 0$. 
\end{thm}

For the special case of (\ref{eq:predictive_density_2_1al_special_case}), we can obtain another sufficient condition. 

\begin{thm}
\label{thm:special_case} 
Suppose that $b = 1$ and $a = p / 2 - 1$ for $p \ge 2$. 
Then $R(( \bmu , \eta ), \ph _{2}^{( \pi _{b, a} )} ) \le R(( \bmu , \eta ), \ph _{2}^{( \pi _0 )} )$ for all $\bmu \in \mathbb{R} ^p$ and $\eta \in (0, \infty )$. 
Equality holds if and only if $p = 2$ and $\bmu = \bm{0} _p$. 
\end{thm}

\section{Simulation Study}
\label{sec:sim}
In this section, we investigate through simulation the numerical performance of the risk functions of the Bayesian predictive densities $\ph _{2}^{\rm{O}} ( \cdot ; \X , V) = \ph _{2}^{( \pi _0 )} ( \cdot ; \X , V)$ and $\ph _{2}^{(b , a)} ( \cdot ; \X , V) = \ph _{2}^{( \pi _{b, a} )} ( \cdot ; \X , V)$ for $b \in \{ n_1 / 2, 1 \} $ and $a \in \{ 0, p / 2 - 1 \} $. 
We consider the following cases: (i) $( n_1 , n_2 ) = (3, 3)$; (ii) $( n_1 , n_2 ) = (3, 5)$; (iii) $( n_1 , n_2 ) = (5, 3)$; (iv) $( n_1 , n_2 ) = (5, 5)$. 
We set $p = 14$. 
When $b = 1$, the conditions of Theorem \ref{thm:KK} are satisfied for both $a = 0$ and $a = p / 2 - 1$. 
On the other hand, when $b = n_1 / 2$, the condition of part (i) of Corollary \ref{cor:D0} is satisfied if $a = p / 2 - 1$ but not if $a = 0$, which can be verified numerically. 

The risk function of $\ph _{2}^{\rm{O}} ( \cdot ; \X , V)$ is a constant independent of the unknown parameters $( \bmu , \eta )$ while that of $\ph _{2}^{(b, a)} ( \cdot ; \X , V)$ depends on $( \bmu , \eta )$ only through $\th = \eta || \bmu ||^2$. 
For $\th \in \{ 0, 20, 40, 60 \} $, we obtain approximated values of the risk function of $\ph _{2}^{(b, a)} ( \cdot ; \X , V)$ by the Monte Carlo simulation with $100,000$ replications. 
The integrals are calculated via the Monte Carlo simulation with $10,000$ replications. 

\begin{figure}
\centering
\includegraphics[width = 16cm]{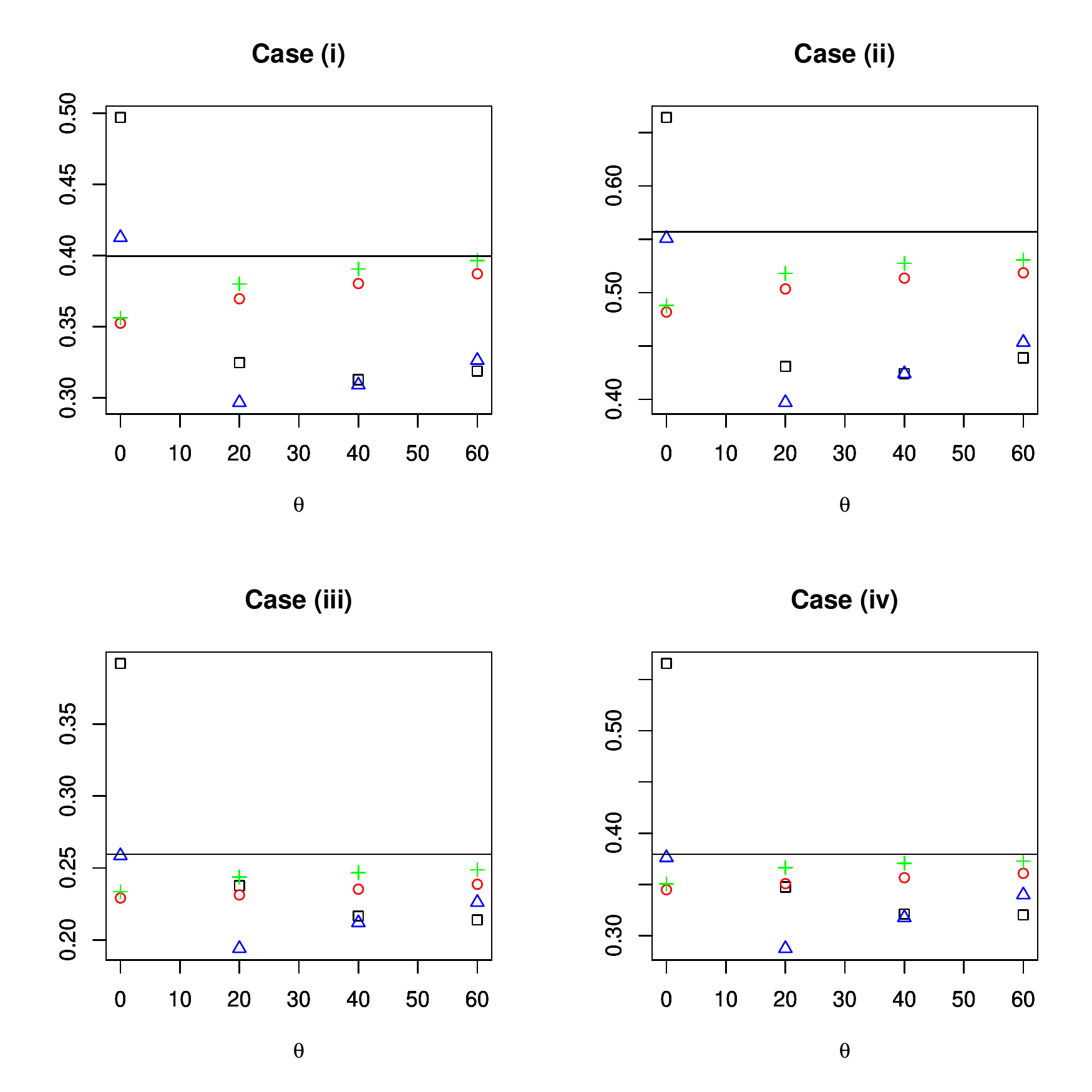}
\caption{Risks of the predictive densities $\ph _{2}^{\rm{O}} ( \cdot ; \X , V)$ and $\ph _{2}^{(b, a)} ( \cdot ; \X , V)$ %
in the following cases: (i) $( n_1 , n_2 ) = (3, 3)$; (ii) $( n_1 , n_2 ) = (3, 5)$; (iii) $( n_1 , n_2 ) = (5, 3)$; (iv) $( n_1 , n_2 ) = (5, 5)$. 
We set $p = 14$. 
The horizontal lines show the constant risk of $\ph _{2}^{\rm{O}} ( \cdot ; \X , V)$. 
The black squares, red circles, blue triangles, and green pluses correspond to $(b, a) = ( n_1 / 2, 0), ( n_1 / 2, p / 2 - 1), (1, 0), (1, p / 2 - 1)$, respectively. }
\label{fig:risk}
\end{figure}

The results are illustrated in Figure \ref{fig:risk}. 
The constant risk of $\ph _{2}^{\rm{O}} ( \cdot ; \X , V)$ is not the same for each case. 
For each $b \in \{ n_1 / 2, 1 \} $, the risk values of $\ph _{2}^{(b, p / 2 - 1)} ( \cdot ; \X , V)$ are smaller than those of $\ph _{2}^{(b, 0)} ( \cdot ; \X , V)$ when $\th = 0$ but larger when $\th = 60$. 
The risk values of $\ph _{2}^{( n_1 / 2, 0)} ( \cdot ; \X , V)$ are larger than those of $\ph _{2}^{(1, 0)} ( \cdot ; \X , V)$ when $\th = 0$ but smaller when $\th = 60$; on the other hand, the risk values of $\ph _{2}^{( n_1 / 2, p / 2 - 1)} ( \cdot ; \X , V)$ are close to those of $\ph _{2}^{(1, p / 2 - 1)} ( \cdot ; \X , V)$ for all $\th \in \{ 0, 20, 40, 60 \} $. 
Since by Theorem \ref{thm:KK} the values of the risk functions of $\ph _{2}^{\rm{O}} ( \cdot ; \X , V)$ and $\ph _{2}^{(1, 0)} ( \cdot ; \X , V)$ at $\th = 0$ coincide, that the blue triangles are not on the horizontal lines when $\th = 0$ will be due to Monte Carlo error. 
Finally, $\ph _{2}^{( n_1 / 2, 0)} ( \cdot ; \X , V)$ does not seem to dominate $\ph _{2}^{\rm{O}} ( \cdot ; \X , V)$ with the value of $a$ too small, for the black squares lie far above the horizontal lines when $\th = 0$.

\section{Appendix}
Useful lemmas are given in Section \ref{subsec:lemmas}. 
Propositions \ref{prp:predictive_density_1} and \ref{prp:predictive_density_2}, Theorems \ref{thm:BM}, \ref{thm:KK}, and \ref{thm:special_case}, and Corollaries \ref{cor:D0} and \ref{cor:2n_2} are proved in Section \ref{subsec:proofs}. 
Let $\mathbb{N} _0 = \{ 0, 1, 2, \dotsc \} $.

\subsection{Lemmas}
\label{subsec:lemmas} 

\begin{lem}
\label{lem:incomplete_beta} 
For any $\xi _1 , \xi _2 , c > 0$, it holds that 
\begin{align}
\int_{0}^{1} {\ga ^{\xi _1 - 1} \over (1 + c \ga )^{\xi _1 + \xi _2}} d\ga = {1 \over c^{\xi _1}} \int_{0}^{c / (1 + c)} \ga ^{\xi _1 - 1} (1 - \ga )^{\xi _2 - 1} d\ga \text{.} \non 
\end{align}
\end{lem}

\noindent
{\bf Proof%
.} \ \ We have 
\begin{align}
\int_{0}^{1} {\ga ^{\xi _1 - 1} \over (1 + c \ga )^{\xi _1 + \xi _2}} d\ga = \int_{0}^{c} {\la ^{\xi _1 - 1} / c^{\xi _1} \over (1 + \la )^{\xi _1 + \xi _2}} d\la = {1 \over c^{\xi _1}} \int_{0}^{c / (1 + c)} \ga ^{\xi _1 - 1} (1 - \ga )^{\xi _2 - 1} d\ga \text{,} \non 
\end{align}
which is the desired result. 
\hfill$\Box$

\begin{lem}
\label{lem:BM_inequality} 
For any $\xi _1 , \xi _{2, 1} , \xi _{2, 2} , c > 0$, we have 
\begin{align}
&\int_{0}^{1} (1 - \ga )^{\xi _{2, 1} - 1} \ga ^{\xi _1 - 1} \Big( 1 - {c \over 1 + c} \ga \Big) ^{\xi _{2, 2}} d\ga \non \\
&\ge B( \xi _{2, 1} + \xi _{2, 2} , \xi _1 ) \times \begin{cases} \displaystyle 1 + \Big\{ {\Ga ( \xi _{2, 1} ) \Ga ( \xi _{2, 1} + \xi _{2, 2} + \xi _1 ) \over \Ga ( \xi _{2, 1} + \xi _{2, 2} ) \Ga ( \xi _{2, 1} + \xi _1 )} - 1 \Big\} {1 \over 1 + c} \text{,} & \text{if $\xi _{2, 2} \le 1$} \text{,} \\ \displaystyle \Big(1 + {\xi _1 \over \xi _{2, 1} + \xi _{2, 2} - 1} {1 \over 1 + c} \Big) ^{\xi _{2, 2}} \text{,} & \text{if $\xi _{2, 2} > 1$} \text{.} \end{cases} \non 
\end{align}
\end{lem}

\noindent
{\bf Proof
.} \ \ Suppose first that $\xi _{2, 2} \le 1$. 
Then, by Lemma 3 of Boisbunon and Maruyama (2014), we have for all $\ga \in (0, 1)$ 
\begin{align}
\Big( 1 - {c \over 1 + c} \ga \Big) ^{\xi _{2, 2}} \ge (1 - \ga )^{\xi _{2, 2}} + \Big( 1 - {c \over 1 + c} \Big) \{ 1 - (1 - \ga )^{\xi _{2, 2}} \} \text{.} \non 
\end{align}
Therefore, 
\begin{align}
&\int_{0}^{1} (1 - \ga )^{\xi _{2, 1} - 1} \ga ^{\xi _1 - 1} \Big( 1 - {c \over 1 + c} \ga \Big) ^{\xi _{2, 2}} d\ga \non \\
&\ge B( \xi _{2, 1} + \xi _{2, 2} , \xi _1 ) + {1 \over 1 + c} \{ B( \xi _{2, 1} , \xi _1 ) - B( \xi _{2, 1} + \xi _{2, 2} , \xi _1 ) \} \non \\
&= B( \xi _{2, 1} + \xi _{2, 2} , \xi _1 ) \Big[ 1 + \Big\{ {\Ga ( \xi _{2, 1}) \Ga ( \xi _{2, 1} + \xi _{2, 2} + \xi _1 ) \over \Ga ( \xi _{2, 1} + \xi _{2, 2} ) \Ga ( \xi _{2, 1} + \xi _1 )} - 1 \Big\} {1 \over 1 + c} \Big] \text{.} \non 
\end{align}
Next suppose that $\xi _{2, 2} > 1$. 
Then, by Jensen's inequality, it follows that 
\begin{align}
&\int_{0}^{1} (1 - \ga )^{\xi _{2, 1} - 1} \ga ^{\xi _1 - 1} \Big( 1 - {c \over 1 + c} \ga \Big) ^{\xi _{2, 2}} d\ga \non \\
&= B( \xi _{2, 1} + \xi _{2, 2} , \xi _1 ) \int_{0}^{1} {(1 - \ga )^{\xi _{2, 1} - 1} \ga ^{\xi _1 - 1} \over B( \xi _{2, 1} + \xi _{2, 2} , \xi _1 )} \Big( 1 - \ga + \ga - {c \over 1 + c} \ga \Big) ^{\xi _{2, 2}} d\ga \non \\
&= B( \xi _{2, 1} + \xi _{2, 2} , \xi _1 ) \int_{0}^{1} {(1 - \ga )^{\xi _{2, 1} + \xi _{2, 2} - 1} \ga ^{\xi _1 - 1} \over B( \xi _{2, 1} + \xi _{2, 2} , \xi _1 )} \Big( 1 + {1 \over 1 + c} {\ga \over 1 - \ga } \Big) ^{\xi _{2, 2}} d\ga \non \\
&\ge B( \xi _{2, 1} + \xi _{2, 2} , \xi _1 ) \Big( 1 + {1 \over 1 + c} {\xi _1 \over \xi _{2, 1} + \xi _{2, 2} - 1} \Big) ^{\xi _{2, 2}} \non 
\end{align}
This completes the proof. 
\hfill$\Box$

\begin{lem}
\label{lem:hypergeo} 
For any $\xi _1 , \xi _2 , c > 0$, we have 
\begin{align}
\int_{0}^{1} \{ \log (1 + c \rho ) \} {\rho ^{\xi _1 - 1} (1 - \rho )^{\xi _2 - 1} \over B( \xi _1 , \xi _2 )} d\rho = \int_{0}^{1} {(1 - \rho )^{\xi _1 + \xi _2 - 1} \over \rho } \Big\{ 1 - {1 \over (1 + c \rho )^{\xi _1}} \Big\} d\rho \text{.} \non 
\end{align}
\end{lem}

\noindent
{\bf Proof%
.} \ \ The hypergeometric function $F$ satisfies 
\begin{align}
F( a' , b' ; c' ; z' ) &= \sum_{s = 0}^{\infty } {\Ga ( a' + s) \Ga ( b' + s) \Ga ( c' ) \over \Ga ( a' ) \Ga ( b' ) \Ga ( c' + s) s!} ( z' )^{s} \non \\
&= {\Ga ( c' ) \over \Ga ( b' ) \Ga ( c' - b' )} \int_{0}^{1} {t^{b' - 1} (1 - t)^{c' - b' - 1} \over (1 - z' t)^{a'}} dt \non 
\end{align}
for $a' > 0$, $c' > b' > 0$, and $z' < 0$. 
Therefore, 
\begin{align}
&\int_{0}^{1} \{ \log (1 + c \rho ) \} {\rho ^{\xi _1 - 1} (1 - \rho )^{\xi _2 - 1} \over B( \xi _1 , \xi _2 )} d\rho - \int_{0}^{1} {(1 - \rho )^{\xi_1 + \xi _2 - 1} \over \rho } \Big\{ 1 - {1 \over (1 + c \rho )^{\xi _1 }} \Big\} d\rho \non \\
&= \int_{0}^{1} \Big( \int_{0}^{1} {c \rho \over 1 + c \rho t} dt \Big) {\rho ^{\xi _1 -1} (1 - \rho )^{\xi _2 - 1} \over B( \xi _1 , \xi _2 )} d\rho - \int_{0}^{1} {(1 - \rho )^{\xi _1 + \xi _2 - 1} \over \rho } \Big\{ \int_{0}^{1} {\xi _1 c \rho \over (1 + c \rho t)^{\xi _1 + 1}} dt \Big\} d\rho \non \\
&= \int_{0}^{1} \Big\{ {1 \over B( \xi _1 , \xi _2 )} \int_{0}^{1} {c \rho ^{\xi _1 } (1 - \rho )^{\xi _2 - 1} \over 1 + c t \rho } d\rho - \xi _1 \int_{0}^{1} {c (1 - \rho )^{\xi _1 + \xi _2 - 1} \over  (1 + c t \rho )^{\xi _1 + 1}} d\rho \Big\} dt \non \\
&= {\xi _1 c \over \xi _1 + \xi _2 } \int_{0}^{1} \{ F(1, \xi _1 + 1; \xi _1 + \xi _2 + 1; - c t) - F( \xi _1 + 1, 1; \xi _1 + \xi _2 + 1; - c t) \} dt = 0 \text{,} \non 
\end{align}
which proves Lemma \ref{lem:hypergeo}. 
\hfill$\Box$

\begin{lem}
\label{lem:digamma} 
For any $\xi _1 , \xi _2 > 0$, we have 
\begin{align}
\psi ( \xi _1 ) - \psi ( \xi _2 ) &= \sum_{i = 0}^{\infty } {\xi _1 - \xi _2 \over (i + \xi _1 ) (i + \xi _2 )} \text{.} \non 
\end{align}
\end{lem}

\noindent
{\bf Proof%
.} \ \ Let $C = \lim_{i \to \infty } \sum_{j = 1}^{i} 1 / j - \log i$. 
Then 
\begin{align}
\psi ( \xi _1 ) - \psi ( \xi _2 ) &= \Big\{ - {1 \over \xi _1} - C + \sum_{i = 1}^{\infty } \Big( {1 \over i} - {1 \over i + \xi _1} \Big) \Big\} - \Big\{ - {1 \over \xi _2} - C + \sum_{i = 1}^{\infty } \Big( {1 \over i} - {1 \over i + \xi _2 } \Big) \Big\} \non \\
&= \sum_{i = 0}^{\infty } \Big( {1 \over i + \xi _2 } - {1 \over i + \xi _1} \Big) = \sum_{i = 0}^{\infty } {\xi _1 - \xi _2 \over (i + \xi _1 ) (i + \xi _2 )} \text{,} \non 
\end{align}
which shows Lemma \ref{lem:digamma}. 
\hfill$\Box$

\begin{lem}
\label{lem:monotonicity} 
Let $\xi _1 > 0$ and $1 < \xi _{2, 1} < \xi _{2, 2}$. 
Let, for $i \in \{ 1, 2 \} $, 
\begin{align}
F_i (q) &= \int_{0}^{q} {\ga ^{\xi _1 - 1} (1 - \ga )^{\xi _{2, i} - 1} \over B( \xi _1 , \xi _{2, i} )} d\ga \text{,} \quad q \in (0, 1) \text{.} \non 
\end{align}
\begin{itemize}
\item[{\rm{(i)}}]
${F_2}^{- 1} ( \om ) / {F_1}^{- 1} ( \om )$ is nondecreasing in $\om \in (0, 1)$. 
\item[{\rm{(ii)}}]
There exist $0 < \underline{\om } < \overline{\om } < 1$ such that ${F_2}^{- 1} ( \om ) / {F_1}^{- 1} ( \om )$ is strictly increasing in $\om \in ( \underline{\om } , \overline{\om } )$. 
\end{itemize}
\end{lem}

\noindent
{\bf Proof%
.} \ \ Part (i) follows from Lemma 2 of Kato (2009). 
For part (ii), we need only show that ${F_2}^{- 1} ( \om ) / {F_1}^{- 1} ( \om )$ is not constant in $\om \in (0, 1)$. 
Suppose that there exists $C_0 \in \mathbb{R}$ such that ${F_2}^{- 1} ( \om ) / {F_1}^{- 1} ( \om ) = C_0$ for all $\om \in (0, 1)$. 
Then $C_0 = \lim_{\om \to 1} \{ {F_2}^{- 1} ( \om ) / {F_1}^{- 1} ( \om ) \} = 1$. 
Therefore, we have that ${F_2}^{- 1} = {F_1}^{- 1}$ and hence that $F_2 = F_1$. 
This is a contradiction. 
\hfill$\Box$

\begin{lem}
\label{lem:gamma_multiplication} 
Let $h \in \mathbb{N}$ and $\xi \ge 1$. 
Then for all $\ta > 0$, 
\begin{align}
{\partial \over \partial \ta } {\Ga ((h + 1) \ta ) \Ga ( \ta + \xi ) \over \Ga ( \ta ) \Ga ((h + 1) \ta + \xi )} \begin{cases} = 0 \text{,} & \text{if $\xi = 1$} \text{,} \\ < 0 \text{,} & \text{if $\xi > 1$} \text{.} \end{cases} \non 
\end{align}
\end{lem}

\noindent
{\bf Proof%
.} \ \ By Gauss's multiplication formula, we have 
\begin{align}
&{\Ga ((h + 1) \ta ) \Ga ( \ta + \xi ) \over \Ga ( \ta ) \Ga ((h + 1) \ta + \xi )} \non \\
&= {\Ga ( \ta + \xi ) \over \Ga ( \ta )} \frac{ (2 \pi )^{\{ 1 - (h + 1) \} / 2} (h + 1)^{(h + 1) \ta - 1 / 2} }{ (2 \pi )^{\{ 1 - (h + 1) \} / 2} (h + 1)^{(h + 1) \ta + \xi - 1 / 2} } {\prod_{i = 0}^{h} \Ga ( \ta + i / (h + 1)) \over \prod_{i = 0}^{h} \Ga ( \ta + \xi / (h + 1) + i / (h + 1))} \non \\
&= {1 \over (h + 1)^{\xi }} {\Ga ( \ta + \xi ) \over \Ga ( \ta + ( \xi + h) / (h + 1))} \prod_{i = 1}^{h} {\Ga ( \ta + i / (h + 1)) \over \Ga ( \ta + ( \xi + i - 1) / (h + 1))} \non 
\end{align}
for all $\ta > 0$. 
Therefore, by Lemma \ref{lem:digamma}, 
\begin{align}
&{\partial \over \partial \ta } \log {\Ga ((h + 1) \ta ) \Ga ( \ta + \xi ) \over \Ga ( \ta ) \Ga ((h + 1) \ta + \xi )} \non \\
&= \psi ( \ta + \xi ) - \psi \Big( \ta + {\xi + h \over h + 1} \Big) + \sum_{i = 1}^{h} \Big\{ \psi \Big( \ta + {i \over h + 1} \Big) - \psi \Big( \ta + {\xi + i - 1 \over h + 1} \Big) \Big\} \non \\
&= {( \xi - 1) h \over h + 1} \sum_{j = 0}^{\infty } \Big\{ \frac{ 1 }{ (j + \ta + \xi ) \big( j + \ta + {\xi + h \over h + 1} \big) } - {1 \over h} \sum_{i = 1}^{h} \frac{ 1 }{ \big( j + \ta + {i \over h + 1} \big) \big( j + \ta + {\xi + i - 1 \over h + 1} \big) } \Big\} \non 
\end{align}
for all $\ta > 0$. 
Fix $j \in \mathbb{N} _0$ and $\ta > 0$. 
Then, by Jensen's inequality, 
\begin{align}
&\frac{ 1 }{ (j + \ta + \xi ) \big( j + \ta + {\xi + h \over h + 1} \big) } - {1 \over h} \sum_{i = 1}^{h} \frac{ 1 }{ \big( j + \ta + {i \over h + 1} \big) \big( j + \ta + {\xi + i - 1 \over h + 1} \big) } \non \\
&\le \frac{ 1 }{ (j + \ta + \xi ) \{ j + \ta + ( \xi + h) / (h + 1) \} } - \frac{ 1 }{ (j + \ta + 1 / 2) \{ j + \ta + ( \xi - 1) / (h + 1) + 1 / 2 \} } \non \\
&= \frac{ (- \xi ) (j + \ta ) + (1 / 2) \{ ( \xi - 1) / (h + 1) + 1 / 2 \} - \xi ( \xi + h) / (h + 1) }{ (j + \ta + \xi ) \{ j + \ta + ( \xi + h) / (h + 1) \} (j + \ta + 1 / 2) \{ j + \ta + ( \xi - 1) / (h + 1) + 1 / 2 \} } < 0 \text{.} \non 
\end{align}
This completes the proof. 
\hfill$\Box$

\subsection{Proofs}
\label{subsec:proofs} 

\noindent
{\bf Proof of Proposition \ref{prp:predictive_density_1}.} \ \ Since the joint posterior density of $( \bmu , \eta )$ is proportional to 
\begin{align}
\eta ^{n_1 / 2 + p / 2 - 1} \exp \Big( - {\eta \over 2} V \Big) \exp \Big( - {\eta \over 2} || \X - \bmu ||^2 \Big) \text{,} \non 
\end{align}
the marginal posterior of $\eta $ is proportional to 
\begin{align}
\eta ^{n_1 / 2 + p / 2 - 1} \exp \Big( - {\eta \over 2} V \Big) \int_{\mathbb{R} ^p} \exp \Big( - {\eta \over 2} || \X - \bmu ||^2 \Big) d\bmu = (2 \pi )^{p / 2} \eta ^{n_1 / 2 - 1} \exp \Big( - {\eta \over 2} V \Big) \text{.} \non 
\end{align}
Therefore, the posterior mean of $p_2 (w | \eta )$ is 
\begin{align}
\ph _{2}^{( \pi _0 )} (w | \X , V) &= {(1 / 2)^{n_2 / 2} \over \Ga ( n_2 / 2)} w^{n_2 / 2 - 1} \frac{ \int_{0}^{\infty } \eta ^{( n_1 + n_2 ) / 2 - 1} e^{- \eta (V + w) / 2} d\eta }{ \int_{0}^{\infty } \eta ^{n_1 / 2 - 1} e^{- \eta V / 2} d\eta } \non \\
&= {(1 / 2)^{n_2 / 2} \over \Ga ( n_2 / 2)} w^{n_2 / 2 - 1} \frac{ \Ga (( n_1 + n_2 ) / 2) / \{ (V + w) / 2 \} ^{( n_1 + n_2 ) / 2} }{ \Ga ( n_1 / 2) / (V / 2)^{n_1 / 2} } \text{,} \non 
\end{align}
which is the desired result. 
\hfill$\Box$

\bigskip

\noindent
{\bf Proof of Proposition \ref{prp:predictive_density_2}.} \ \ Let $\pi _{b, a} ( \ga ) = (1 - \ga )^{b - 1} \ga ^{- a - 1}$ for $\ga \in (0, 1)$. 
Then the joint posterior density of $( \bmu , \eta )$ is proportional to 
\begin{align}
&\int_{0}^{1} \pi _{b, a} ( \ga ) \Big( {\ga \over 1 - \ga } \Big) ^{p / 2} \eta ^{n_1 / 2 + p - a - 1} \exp \Big( - {\eta \over 2} V \Big) \exp \Big\{ - {\eta \over 2} \Big( {\ga \over 1 - \ga } || \bmu ||^2 + || \X - \bmu ||^2 \Big) \Big\} d\ga \text{.} \non 
\end{align}
Note that 
\begin{align}
{\ga \over 1 - \ga } || \bmu ||^2 + || \X - \bmu ||^2 &= {|| \bmu - (1 - \ga ) \X ||^2 \over 1 - \ga } + \ga || \X ||^2 \text{.} \non 
\end{align}
Then the marginal posterior of $\eta $ is proportional to 
\begin{align}
&\int_{0}^{1} \pi _{b, a} ( \ga ) \Big( {\ga \over 1 - \ga } \Big) ^{p / 2} \eta ^{n_1 / 2 + p - a - 1} \exp \Big( - {\eta \over 2} V \Big) \Big( \int_{\mathbb{R} ^p} \exp \Big[ - {\eta \over 2} \Big\{ {|| \bmu - (1 - \ga ) \X ||^2 \over 1 - \ga } + \ga || \X ||^2 \Big\} \Big] d\bmu \Big) d\ga \non \\
&= (2 \pi )^{p / 2} \int_{0}^{1} \pi _{b, a} ( \ga ) \ga ^{p / 2} \eta ^{n_1 / 2 + p / 2 - a - 1} \exp \Big\{ - {\eta \over 2} (V + \ga || \X ||^2 ) \Big\} d\ga \text{.} \non 
\end{align}
Therefore, the Bayesian predictive density $\ph _{2}^{( \pi _{b, a} )} ( \cdot | \X , V)$ is given by 
\begin{align}
\frac{ \displaystyle \ph _{2}^{( \pi _{b, a} )} (w | \X , V) }{ \displaystyle {(1 / 2)^{n_2 / 2} \over \Ga ( n_2 / 2)} w^{n_2 / 2 - 1} } &= \frac{ \displaystyle \int_{0}^{1} \pi _{b, a} ( \ga ) \ga ^{p / 2} \Big[ \int_{0}^{\infty } \eta ^{( n_1 + n_2 ) / 2 + p / 2 - a - 1} \exp \Big\{ - {\eta \over 2} (V + w + \ga || \X ||^2 ) \Big\} d\eta \Big] d\ga }{ \displaystyle \int_{0}^{1} \pi _{b, a} ( \ga ) \ga ^{p / 2} \Big[ \int_{0}^{\infty } \eta ^{n_1 / 2 + p / 2 - a - 1} \exp \Big\{ - {\eta \over 2} (V + \ga || \X ||^2 ) \Big\} d\eta \Big] d\ga } \non \\
&= \frac{ \displaystyle \int_{0}^{1} \pi _{b, a} ( \ga ) \ga ^{p / 2} {\Ga (( n_1 + n_2 ) / 2 + p / 2 - a) \over \{ (1 / 2) (V + w + \ga || \X ||^2 ) \} ^{( n_1 + n_2 ) / 2 + p / 2 - a}} d\ga }{ \displaystyle \int_{0}^{1} \pi _{b, a} ( \ga ) \ga ^{p / 2} {\Ga ( n_1 / 2 + p / 2 - a) \over \{ (1 / 2) (V + \ga || \X ||^2 ) \} ^{n_1 / 2 + p / 2 - a}} d\ga } \text{,} \non 
\end{align}
from which the desired result follows. 
\hfill$\Box$

\bigskip

\noindent
{\bf Proof of Theorem \ref{thm:BM}.} \ \ Let $\De = R(( \bmu , \eta ), \ph _{2}^{( \pi _{n_1 / 2, a} )} ) - R(( \bmu , \eta ), \ph _{2}^{( \pi _0 )} )$. 
By Propositions \ref{prp:predictive_density_1} and \ref{prp:predictive_density_2} and by (\ref{eq:denominator}) and (\ref{eq:numerator}), we have 
\begin{align}
\De &= E_{( \bmu , \eta )}^{( \X , V, W)} \Big[ \log {\ph _{2}^{( \pi _0 )} (W; \X , V) \over \ph _{2}^{( \pi _{n_1 / 2, a} )} (W; \X , V)} \Big] \non \\
&= E_{( \bmu , \eta )}^{( \X , V, W)} \Big[ \log B \Big( {n_1 + n_2 \over 2}, {p \over 2} - a \Big) + \Big( {p \over 2}  - a \Big) \log {V + W + || \X ||^2 \over V + || \X ||^2} \non \\
&\quad - \log \int_{0}^{1} (1 - \ga )^{n_1 / 2 - 1} \ga ^{p / 2 - a - 1} \Big( 1 - {|| \X ||^2 \over V + W + || \X ||^2} \ga \Big) ^{n_2 / 2} d\ga \Big] \text{.} \non 
\end{align}
It follows from Lemma \ref{lem:BM_inequality} that for all $\x \in \mathbb{R} ^p$, $v \in (0, \infty )$, and $w \in (0, \infty )$, 
\begin{align}
&\int_{0}^{1} (1 - \ga )^{n_1 / 2 - 1} \ga ^{p / 2 - a - 1} \Big( 1 - {|| \x ||^2 \over v + w + || \x ||^2} \ga \Big) ^{n_2 / 2} d\ga \non \\
&\ge B \Big( {n_1 + n_2 \over 2}, {p \over 2} - a \Big) \Big( 1 + c_1 {v + w \over v + w + || \x ||^2} \Big) ^{c_2} \text{.} \non 
\end{align}
Therefore, 
\begin{align}
\De &\le E_{( \bmu , \eta )}^{( \X , V, W)} \Big[ \Big( {p \over 2}  - a \Big) \log {V + W + || \X ||^2 \over V + || \X ||^2} - c_2 \log \Big( 1 + c_1 {V + W \over V + W + || \X ||^2} \Big) \Big] \non \\
&= E_{( \bmu , \eta )}^{( \X , V, W)} \Big[ \Big( {p \over 2}  - a \Big) \log {\eta V + \eta W + || \sqrt{\eta } \X ||^2 \over \eta V + || \sqrt{\eta } \X ||^2} - c_2 \log \Big( 1 + c_1 {\eta V + \eta W \over \eta V + \eta W + || \sqrt{\eta } \X ||^2} \Big) \Big] \text{.} \non 
\end{align}
Let $k = n_1 / 2$, $l = n_2 / 2$, $m = p / 2$, and $m' = m - a = p / 2 - a$. 
Let $Z \sim {\rm{Po}} ( \th / 2)$ for $\th = \eta || \bmu ||^2$ and let $\Vt $, $\Wt $, and $\Tt $ be independently distributed as $\chi ^2 ( n_1 )$, $\chi ^2 ( n_2 )$, and $\chi ^2 (p + 2 Z)$, respectively. 
Then since $( \eta V, \eta W, || \sqrt{\eta } \X ||^2 ) \stackrel{\rm{d}}{=} ( \Vt , \Wt , \Tt )$ and since the expectation of the logarithm of a Chi-squared variable with $\nu > 0$ degrees
of freedom is $\log 2 + \psi ( \nu / 2)$, it follows that 
\begin{align}
\De &\le E_{\th }^{Z} \Big[ E_{\th }^{( \Tt , \Vt , \Wt ) | Z} \Big[ m' \log {\Vt + \Wt + \Tt \over \Vt + \Tt } - c_2 \log \Big( 1 + c_1 {\Vt + \Wt \over \Vt + \Wt + \Tt } \Big) \Big| Z \Big] \Big] \non \\
&= E_{\th }^{Z} [ D_1 (Z) + D_2 (Z) ] \text{,} \label{tBMp1} 
\end{align}
where 
\begin{align}
D_1 (z) &= m' \{ \psi (k + l + m + z) - \psi (k + m + z) \} \text{,} \quad z \in \mathbb{N} _0 \text{,} \non 
\end{align}
and 
\begin{align}
D_2 (z) &= E_{\th }^{\rho _Z | Z} [ - c_2 \log (1 + c_1 \rho _Z ) | Z = z ] \text{,} \quad z \in \mathbb{N} _0 \text{,} \non 
\end{align}
for a random variable $\rho _Z$ such that $\rho _Z | Z \sim {\rm{Beta}} (k + l, m + Z)$. 
By Lemma \ref{lem:hypergeo}, 
\begin{align}
D_2 (z) &= - c_2 \int_{0}^{1} \{ \log (1 + c_1 \rho ) \} {\rho ^{k + l - 1} (1 - \rho )^{m + z - 1} \over B(k + l, m + z)} d\rho \non \\
&= - c_2 \int_{0}^{1} {(1 - \rho )^{k + l + m + z - 1} \over \rho } \Big\{ 1 - {1 \over (1 + c_1 \rho )^{k + l}} \Big\} d\rho \non 
\end{align}
for all $z \in \mathbb{N} _0$. 
Therefore, by Lemma \ref{lem:digamma}, $\lim_{z \to \infty } \{ D_1 (z) + D_2 (z) \} = 0$. 
Fix $z \in \mathbb{N} _0$. 
Then 
\begin{align}
&\{ D_1 (z + 1) + D_2 (z + 1) \} - \{ D_1 (z) + D_2 (z) \} \non \\
&= m' \Big( {1 \over k + l + m + z} - {1 \over k + m + z} \Big) - c_2 \int_{0}^{1} {(1 - \rho )^{k + l + m + z - 1} (- \rho ) \over \rho } \Big\{ 1 - {1 \over (1 + c_1 \rho )^{k + l}} \Big\} d\rho \non \\
&= - {l m' \over (k + m + z) (k + l + m + z)} + c_2 \int_{0}^{1} (1 - \rho )^{k + l + m + z - 1} \Big\{ 1 - {1 \over (1 + c_1 \rho )^{k + l}} \Big\} d\rho \text{.} \non 
\end{align}
Therefore, 
\begin{align}
D_1 (z + 1) + D_2 (z + 1) \gtreqless D_1 (z) + D_2 (z) \quad \text{if and only if} \quad f(k + m + z) \gtreqless l m' / c_2 \non 
\end{align}
for the function $f$ defined by 
\begin{align}
f( \ze ) &= \ze ( \ze + l) \int_{0}^{1} (1 - \rho )^{\ze + l - 1} \Big\{ 1 - {1 \over (1 + c_1 \rho )^{k + l}} \Big\} d\rho \text{,} \quad \ze \in (0, \infty ) \text{.} \non 
\end{align}
Furthermore, by integration by parts 
\begin{align}
f( \ze ) &= \ze \int_{0}^{1} (1 - \rho )^{\ze + l} {\partial \over \partial \rho } \Big\{ 1 - {1 \over (1 + c_1 \rho )^{k + l}} \Big\} d\rho = \int_{0}^{1} \ze (1 - \rho )^{\ze - 1} (1 - \rho )^{l + 1} {(k + l) c_1 \over (1 + c_1 \rho )^{k + l + 1}} d\rho \non \\
&= \Big[ - (1 - \rho )^{\ze } (1 - \rho )^{l + 1} {(k + l) c_1 \over (1 + c_1 \rho )^{k + l + 1}} \Big] _{0}^{1} + \int_{0}^{1} (1 - \rho )^{\ze } {\partial \over \partial \rho } \Big\{ (1 - \rho )^{l + 1} {(k + l) c_1 \over (1 + c_1 \rho )^{k + l + 1}} \Big\} d\rho \non 
\end{align}
for all $\ze \in (0, \infty )$ and thus $f$ is an increasing function. 
Finally, $D_1 (0) + D_2 (0) \le 0$ by assumption. 
Hence, we conclude that $D_1 (z) + D_2 (z) \le 0$ for all $z \in \mathbb{N} _0$ with strict inequality for some $z \in \mathbb{N} _0$. 
This completes the proof. 
\hfill$\Box$

\bigskip

\noindent
{\bf Proof of Corollary \ref{cor:D0}.} \ \ Let $k = n_1 / 2$, $l = n_2 / 2$, $m = p / 2$, $m' = p / 2 - a$. 
We show that 
\begin{align}
{m' \over c_2} \{ \psi (k + l + m) - \psi (k + m) \} \le \int_{0}^{1} (1 - \rho )^{k + l + m - 1} g( \rho ) d\rho \text{,} \label{cD0p1} 
\end{align}
where 
\begin{align}
( c_1 , c_2 ) &= \begin{cases} \displaystyle \Big( {\Ga (k) \Ga (k + l + m' ) \over \Ga (k + l) \Ga (k + m' )} - 1, 1 \Big) \text{,} & \text{if $l \le 1$} \text{,} \\ \displaystyle \Big( {m' \over k + l - 1}, l \Big) \text{,} & \text{if $l > 1$} \text{,} \end{cases} \non 
\end{align}
and where $g \colon (0, 1) \to [0, \infty )$ is the function defined by 
\begin{align}
g( \rho ) &= {1 \over \rho } \Big\{ 1 - {1 \over (1 + c_1 \rho )^{k + l}} \Big\} = {1 \over \rho } \int_{0}^{1} \Big[ {\partial \over \partial t} \Big\{ {- 1 \over (1 + c_1 \rho t)^{k + l}} \Big\} \Big] dt = \int_{0}^{1} {(k + l) c_1 \over (1 + c_1 \rho t)^{k + l + 1}} dt \text{,} \quad \rho \in (0, 1) \text{.} \label{cD0p2} 
\end{align}

For part (i), since for all $\rho \in (0, 1)$ 
\begin{align}
g' ( \rho ) &= \int_{0}^{1} {- (k + l + 1) (k + l) {c_1}^2 t \over (1 + c_1 \rho t)^{k + l + 2}} dt \text{,} \non 
\end{align}
$g$ is a convex function. 
Therefore, by Jensen's inequality, 
\begin{align}
\int_{0}^{1} (1 - \rho )^{k + l + m - 1} g( \rho ) d\rho &= B(1, k + l + m) \int_{0}^{1} {\rho ^{1 - 1} (1 - \rho )^{k + l + m - 1} \over B(1, k + l + m)} g( \rho ) d\rho \non \\
&\ge B(1, k + l + m) g \Big( {1 \over k + l + m + 1} \Big) \non \\
&= {1 \over k + l + m} \int_{0}^{1} {(k + l) c_1 \over [1 + \{ c_1 / (k + l + m + 1) \} t]^{k + l + 1}} dt \non \\
&= {k + l + m + 1 \over k + l + m} \Big[ 1 - {1 \over \{ 1 + c_1 / (k + l + m + 1) \} ^{k + l}} \Big] \text{,} \non 
\end{align}
the right-hand side of which is greater than or equal to the left-hand side of (\ref{cD0p1}) by assumption. 

To prove part (ii), note that 
\begin{align}
\lim_{m' \to 0} {g( \rho ) \over m'} &= \lim_{m' \to 0} {c_1 \over m'} \int_{0}^{1} {k + l \over (1 + c_1 \rho t)^{k + l + 1}} dt = (k + l) \lim_{m' \to 0} {c_1 \over m'} \non \\
&= (k + l) \begin{cases} \displaystyle {\partial \over \partial m'} \Big| _{m' = 0} {\Ga (k) \Ga (k + l + m' ) \over \Ga (k + l) \Ga (k + m' )} \text{,} & \text{if $l \le 1$} \text{,} \\ \displaystyle {1 \over k + l - 1} \text{,} & \text{if $l > 1$} \text{,} \end{cases} \non \\
&= (k + l) \begin{cases} \displaystyle \psi (k + l) - \psi (k) \text{,} & \text{if $l \le 1$} \text{,} \\ \displaystyle {1 \over k + l - 1} \text{,} & \text{if $l > 1$} \text{,} \end{cases} \non 
\end{align}
Then 
\begin{align}
&\lim_{m' \to 0} {c_2 \over m'} \int_{0}^{1} (1 - \rho )^{k + l + m - 1} g( \rho ) d\rho \non \\
&= c_2 \Big\{ \int_{0}^{1} (1 - \rho )^{k + l + m - 1} d\rho \Big\} (k + l) \begin{cases} \displaystyle \psi (k + l) - \psi (k) \text{,} & \text{if $l \le 1$} \text{,} \\ \displaystyle {1 \over k + l - 1} \text{,} & \text{if $l > 1$} \text{,} \end{cases} \non \\
&= {(k + l) c_2 \over k + l + m} \begin{cases} \displaystyle \psi (k + l) - \psi (k) \text{,} & \text{if $l \le 1$} \text{,} \\ \displaystyle {1 \over k + l - 1} \text{,} & \text{if $l > 1$} \text{,} \end{cases} \non 
\end{align}
from which the desired result follows. 
\hfill$\Box$

\bigskip

\noindent
{\bf Proof of Corollary \ref{cor:2n_2}.} \ \ Let $\De $ and $D_1 (z), D_2 (z)$, $z \in \mathbb{N} _0$, be defined as in the proof of Theorem \ref{thm:BM}. 
For part (i), note that equality holds in (\ref{tBMp1}) when $n_2 = 2$. 
Then if $\ph _{2}^{( \pi _{n_1 / 2, a} )}$ dominates $\ph _{2}^{( \pi _0 )}$, we have $\De |_{\bmu = \bm{0} _p} \le 0$, which implies $D_1 (0) + D_2 (0) \le 0$. 
This proves the ``only if'' part. 
The ``if'' part follows from Theorem \ref{thm:BM}. 
For part (ii), note that by (\ref{cD0p2}), the right-hand side of (\ref{eq:c2n_2-1}) divided by $p / 2 - a$ is 
\begin{align}
&\int_{0}^{1} (1 - \rho )^{n_1 / 2 + p / 2} \Big( \int_{0}^{1} {( n_1 + 2) / n_1 \over [1 + \{ (p - 2 a) / n_1 \} \rho t]^{n_1 / 2 + 2}} dt \Big) d\rho \text{.} \non 
\end{align}
Since the above integral is increasing in $a$, we need only show that equality holds in (\ref{eq:c2n_2-1}) when $a = 0$. 
Suppose that $n_1 = 2$ and that $a = 0$. 
Let $m = p / 2$. 
Then, by integration by parts, 
\begin{align}
&{1 \over p / 2 - a} \int_{0}^{1} (1 - \rho )^{n_1 / 2 + p / 2} {1 \over \rho } \Big( 1 - {1 \over [1 + \{ (p / 2 - a) / ( n_1 / 2) \} \rho ]^{n_1 / 2 + 1}} \Big) d\rho \Big\} \non \\
&= \int_{0}^{1} {(1 - \rho )^{m + 1} \over m \rho } \Big\{ 1 - {1 \over (1 + m \rho )^2} \Big\} d\rho = \int_{0}^{1} (1 - \rho )^{m + 1} \Big\{ {1 \over 1 + m \rho } + {1 \over (1 + m \rho )^2} \Big\} d\rho \non \\
&= \int_{0}^{1} {(1 - \rho )^{m + 1}\over 1 + m \rho } d\rho + {1 \over m} - {m + 1 \over m} \int_{0}^{1} {(1 - \rho )^{m} \over 1 + m \rho } d\rho = {1 \over m} - {1 \over m} \int_{0}^{1} (1 - \rho )^m d\rho = {1 \over m + 1} \text{,} \non 
\end{align}
which equals $1 / ( n_1 / 2 + p / 2)$. 
Thus, we have proved the desired result. 
\hfill$\Box$

\bigskip

\noindent
{\bf Proof of Theorem \ref{thm:KK}.} \ \ By %
(\ref{eq:predictive_density_2_1al}), we have 
\begin{align}
R(( \bmu , \eta ), \ph _{2}^{( \pi _{1, a} )} ) - R(( \bmu , \eta ), \ph _{2}^{( \pi _0 )} ) &= E_{( \bmu , \eta )}^{( \X , V, W)} \Big[ \log {\ph _{2}^{( \pi _0 )} (W; \X , V) \over \ph _{2}^{( \pi _{1, a} )} (W; \X , V)} \Big] = \De ( n_1 + n_2 ) - \De ( n_1 ) \text{,} \non 
\end{align}
where, for each $n \in \{ n_1 , n_1 + n_2 \} $, 
\begin{align}
\De (n) &= E_{( \bmu , \eta )}^{( \X , U_n )} \Big[ - \log \int_{0}^{|| \X ||^2 / ( U_n + || \X ||^2 )} {\ga ^{p / 2 - a - 1} (1 - \ga )^{n / 2 - 1} \over B(p / 2 - a, n / 2)} d\ga \Big] \non 
\end{align}
for the random variable $U_n$ which is $V$ if $n = n_1$ and $V + W$ if $n = n_1 + n_2$. 
Let $\th $, $Z$, and $\Tt $ be defined as in the proof of Theorem \ref{thm:BM}. 
Then for $n \in \{ n_1 , n_1 + n_2 \} $, $\De (n)$ can be written as 
\begin{align}
\De (n) &= E_{\th }^{Z} [ D(n; Z) ] \text{,} \non 
\end{align}
where 
\begin{align}
D(n; z) &= E_{\th }^{( \Tt , \Ut _n ) | Z} \Big[ - \log \int_{0}^{\Tt / ( \Ut _n + \Tt )} {\ga ^{p / 2 - a - 1} (1 - \ga )^{n / 2 - 1} \over B(p / 2 - a, n / 2)} d\ga \Big| Z = z \Big] \text{,} \quad z \in \mathbb{N} _0 \text{,} \non 
\end{align}
for an independent variable $\Ut _n \sim \chi ^2 (n)$. 

Fix $z \in \mathbb{N} _0$. 
Then for each $n \in \{ n_1 , n_1 + n_2 \} $, since $\{ \Tt / ( \Ut _n + \Tt ) \} | (Z = z) \sim {\rm{Beta}} (p / 2 + z, n / 2)$, it follows that 
\begin{align}
D(n; z) %
&= - \int_{0}^{1} \Big\{ \log \int_{0}^{q} {\ga ^{p / 2 - a - 1} (1 - \ga )^{n / 2 - 1} \over B(p / 2 - a, n / 2)} d\ga \Big\} {q^{p / 2 + z - 1} (1 - q)^{n / 2 - 1} \over B(p / 2 + z, n / 2)} dq \non \\
&= - \int_{0}^{1} ( \log \om ) {B(p / 2 - a, n / 2) \over B(p / 2 + z, n / 2)} \{ {F_n}^{- 1} ( \om ) \} ^{z + a} d\om \text{,} \non 
\end{align}
where 
\begin{align}
F_n (q) &= \int_{0}^{q} {\ga ^{p / 2 - a - 1} (1 - \ga )^{n / 2 - 1} \over B(p / 2 - a, n / 2)} d\ga \non 
\end{align}
for $q \in (0, 1)$. 
Therefore, $D( n_1 + n_2 ; z) \lesseqgtr D( n_1 ; z)$ if and only if 
\begin{align}
&\int_{0}^{1} ( \log \om ) \Big[ 1 - C(z) \Big\{ {{F_{n_1}}^{- 1} ( \om ) \over {F_{n_1 + n_2}}^{- 1} ( \om )} \Big\} ^{z + a} \Big] d{P_z} ( \om ) \gtreqless 0 \text{,} \label{tKKp1} 
\end{align}
where 
\begin{align}
C(z) &= {B(p / 2 - a, n_1 / 2) \over B(p / 2 + z, n_1 / 2)} / {B(p / 2 - a, ( n_1 + n_2 ) / 2) \over B(p / 2 + z, ( n_1 + n_2 ) / 2)} \non 
\end{align}
and where $P_z$ is the probability measure with density 
\begin{align}
{B(p / 2 - a, ( n_1 + n_2 ) / 2) \over B(p / 2 + z, ( n_1 + n_2 ) / 2)} \{ {F_{n_1 + n_2}}^{- 1} ( \om ) \} ^{z + a} \text{,} \quad \om \in (0, 1) \text{.} \non 
\end{align}
Since $a < p / 2$ and $n_1 > 2$ by assumption, it follows from Lemma \ref{lem:monotonicity} that ${F_{n_1 + n_2}}^{- 1} ( \om ) / {F_{n_1}}^{- 1} ( \om )$ is nondecreasing in $\om \in (0, 1)$ and strictly increasing in $\om \in ( \underline{\om } , \overline{\om } )$ for some $0 < \underline{\om } < \overline{\om } < 1$. 
Thus, since 
\begin{align}
\int_{0}^{1} \Big[ 1 - C(z) \Big\{ {{F_{n_1}}^{- 1} ( \om ) \over {F_{n_1 + n_2}}^{- 1} ( \om )} \Big\} ^{z + a} \Big] dP( \om ) = 0 \text{,} \non 
\end{align}
the left-hand side of (\ref{tKKp1}) is, by the covariance inequality, greater than zero if $z + a > 0$ and equal to zero if $z + a = 0$, from which the desired result follows. 
\hfill$\Box$

\bigskip

\noindent
{\bf Proof of Theorem \ref{thm:special_case}.} \ \ Let $\th $ and $Z$ be defined as in the proof of Theorem \ref{thm:BM}. 
Then, by the proof of Theorem \ref{thm:KK}, 
\begin{align}
&R(( \bmu , \eta ), \ph _{2}^{( \pi _{1, p / 2 - 1} )} ) - R(( \bmu , \eta ), \ph _{2}^{( \pi _0 )} ) \non \\
&= E_{\th }^{Z} \Big[ \Big[ - \int_{0}^{1} \Big\{ \log \int_{0}^{q} {(1 - \ga )^{n / 2 - 1} \over B(1, n / 2)} d\ga \Big\} {q^{p / 2 + Z - 1} (1 - q)^{n / 2 - 1} \over B(p / 2 + Z, n / 2)} dq \Big] _{n = n_1}^{n = n_1 + n_2} \Big] \non \\
&= E_{\th }^{Z} \Big[ \Big[ - \int_{0}^{1} [ \log \{ 1 - (1 - q)^{n / 2} \} ] {q^{p / 2 + Z - 1} (1 - q)^{n / 2 - 1} \over B(p / 2 + Z, n / 2)} dq \Big] _{n = n_1}^{n = n_1 + n_2} \Big] \non \\
&= \sum_{h = 1}^{\infty } {1 \over h} E_{\th }^{Z} \Big[ \Big[ \int_{0}^{1} {q^{p / 2 + Z - 1} (1 - q)^{(h + 1) (n / 2) - 1} \over B(p / 2 + Z, n / 2)} dq \Big] _{n = n_1}^{n = n_1 + n_2} \Big] \text{.} \non 
\end{align}
Therefore, 
\begin{align}
&R(( \bmu , \eta ), \ph _{2}^{( \pi _{1, p / 2 - 1} )} ) - R(( \bmu , \eta ), \ph _{2}^{( \pi _0 )} ) \non \\
&= \sum_{h = 1}^{\infty } {1 \over h} E_{\th }^{Z} \Big[ \Big[ {B(p / 2 + Z, (h + 1) (n / 2)) \over B(p / 2 + Z, n / 2)} \Big] _{n = n_1}^{n = n_1 + n_2} \Big] \non \\
&= \sum_{h = 1}^{\infty } {1 \over h} E_{\th }^{Z} \Big[ \Big[ {\Ga (p / 2 + Z + n / 2) \Ga ((h + 1) (n / 2)) \over \Ga (p / 2 + Z + (h + 1) (n / 2)) \Ga (n / 2)} \Big] _{n = n_1}^{n = n_1 + n_2} \Big] \non \\
&= \sum_{h = 1}^{\infty } {1 \over h} E_{\th }^{Z} \Big[ \int_{n_1 / 2}^{( n_1 + n_2 ) / 2} \Big\{ {\partial \over \partial \ta } {\Ga (p / 2 + Z + \ta ) \Ga ((h + 1) \ta ) \over \Ga (p / 2 + Z + (h + 1) \ta ) \Ga ( \ta )} \Big\} d\ta \Big] \text{.} \non 
\end{align}
Thus, by Lemma \ref{lem:gamma_multiplication}, we have $R(( \bmu , \eta ), \ph _{2}^{( \pi _{1, p / 2 - 1} )} ) \le R(( \bmu , \eta ), \ph _{2}^{( \pi _0 )} )$. 
Equality holds if and only if $p = 2$ and $\bmu = \bm{0} _p$. 
This completes the proof. 
\hfill$\Box$

\section*{Acknowledgments}
%We would like to thank 
Research of the authors was supported in part by Grant-in-Aid for Scientific Research 
(20J10427, 18K11188) from Japan Society for the Promotion of Science.


\begin{thebibliography}{00}


\bibitem{a1975}
Aitchison, J. (1975). 
Goodness of prediction fit. 
{\it Biometrika}, {\bf 62}, 547--554. 

\bibitem{bm2014}
Boisbunon, A. and Maruyama, Y. (2014). 
Inadmissibility of the best equivariant predictive density in the unknown variance case. 
{\it Biometrika}, {\bf 101}, 733--740.

\bibitem{b1974}
Brewster, J.F. and Zidek, J.V. (1974). 
Improving on equivariant estimators. 
{\it Annals of Statistics}, {\bf 2}, 21--38. 

\bibitem{bgx2008}
Brown, L.D., George, E.I., and Xu, X. (2008).
Admissible predictive density estimation.
{\it Annals of Statistics}, {\bf 36}, 1156--1170. 

\bibitem{glx2006}
George, E.I., Liang, F. and Xu, X. (2006). 
Improved minimax predictive densities under Kullback-Leibler loss. 
{\it Annals of Statistics}, {\bf 34}, 78--91.

\bibitem{glx2012}
George, E.I., Liang, F. and Xu, X. (2012). 
From minimax shrinkage estimation to minimax shrinkage prediction. 
{\it Statistical Science}, {\bf 27}, 82--94. 

\bibitem{k2009}
Kato, K. (2009). 
Improved prediction for a multivariate normal distribution with unknown mean and variance. 
{\it Annals of the Institute of Statistical Mathematics}, {\bf 61}, 531--542. 

\bibitem{k2001}
Komaki, F. (2001). 
A shrinkage predictive distribution for multivariate normal observables. 
{\it Biometrika}, {\bf 88}, 859--864. 

\bibitem{k2009}
Komaki, F. (2009). 
Bayesian predictive densities based on superharmonic priors for the 2-dimensional Wishart model. 
{\it Journal of Multivariate Analysis}, {\bf 100}, 2137--2154. 

\bibitem{k1994}
Kubokawa, T. (1994). 
A unified approach to improving equivariant estimators. 
{\it Annals of Statistics}, {\bf 22}, 290--299. 

\bibitem{lb2004}
Liang, F. and Barron, A. (2004). 
Exact minimax strategies for predictive density estimation, data compression, and model selection. 
{\it IEEE Transactions on Information Theory}, {\bf 50}, 2708--2726. 

\bibitem{l2017}
L'Moudden, A., Marchand, \'E., Kortbi, O., and Strawderman, W.E. (2017). 
On Predictive density estimation for Gamma models with parametric constraints. 
{\it Journal of Statistical Planning and Inference}, {\bf 185}, 56--68. 

\bibitem{m1998}
Maruyama, Y. (1998). 
Minimax estimators of a normal variance. 
{\it Metrika}, {\bf 48}, 209--214. 

\bibitem{ms2005}
Maruyama, Y. and Strawderman W.E. (2005). 
A new class of generalized Bayes minimax ridge regression estimators. 
{\it Annals of Statistics}, {\bf 33}, 1753--1770. 

\bibitem{ms2012}
Maruyama, Y. and Strawderman W.E. (2012). 
Bayesian predictive densities for linear regression models under $\al $-divergence loss: Some results and open problems. 
In {\it IMS Collections, Contemporary Developments in Bayesian
Analysis and Statistical Decision Theory: A Festschrift for William E. Strawderman}, D. Fourdrinier, \'{E}. Marchand \& A. Rukhin, eds., vol. 8. Beachwood, USA: Institute of Mathematical Statistics, 42--56. 

\bibitem{ms2020a}
Maruyama, Y. and Strawderman W.E. (2020a). 
Admissible Bayes equivariant estimation of location vectors for spherically symmetric distributions with unknown scale. 
{\it Annals of Statistics}, {\bf 48}, 1052--1071. 

\bibitem{ms2020b}
Maruyama, Y. and Strawderman W.E. (2020b). 
Admissible estimators of a multivariate normal mean vector when the scale is unknown. 
{\it arXiv preprint arXiv:2003.08571}. 

\bibitem{s1964}
Stein, C. (1964). 
Inadmissibility of the usual estimator for the variance of a normal distribution with unknown mean. 
{\it Annals of the Institute of Statistical Mathematics}, {\bf 16}, 155--160. 




\end{thebibliography}
\end{document}